\documentclass[11pt]{article}
 \usepackage{amsmath,amssymb,amscd}
\newcommand{\pf}[1]{\trivlist \item[\hskip\labelsep\it #1\ ]}
\newcommand{\varpf}[1]{\trivlist \item[\hskip\labelsep\sc #1:]}
\newcommand{\qedbox}{$\rlap{$\sqcap$}\sqcup$}
\newcommand{\qed}{\qquad \qedbox \endtrivlist}
\newcommand{\varqed}{\hfill \rule{0.6em}{0.6em} \endtrivlist}

\newtheorem{theorem}{Theorem}[section]

\newtheorem{df}[theorem]{Definition}

\begin{document}
\title{Weierstra{\ss}\footnote{The last letter is not a $\beta$, but an ``{\ss}'', a letter unique to the German alphabet, which is pronounced like an ``s''.}}
\author{\it Volker Runde}
\date{}
\maketitle
Each university and each department, I guess, develops a peculiar kind of folklore: anecdotes about those of its graduates (or dropouts)
that somehow managed to become famous (or notorious). Very often, there is an element of glee to these stories: ``Well, he may now be a government minister, but {\it I\/} flunked 
him in calculus!'' And also very often, it is impossible to tell truth from legend.
\par
When I was a math student at M\"unster, Germany, in the 1980s, such anecdotes centered mainly around two people: Gerd Faltings, the first and only German to win the Fields medal, mathematics'
equivalent of the Nobel prize, and Karl Weierstra{\ss}\footnote{Books written in English usually spell the name ``Weierstrass''.}, the man who (besides many other mathematical
accomplishments) introduced $\epsilon$s and $\delta$s into calculus. Weierstra{\ss} had been a student at M\"unster in the 1830s and 1840s. There was no one around anymore
who knew anybody who had known anybody who had known anybody who had known Weierstra{\ss}. But this didn't prevent the folklore from blooming. According to folklore, Weierstra{\ss}
flunked out of law school (because he spent most of his time there drinking beer and doing mathematics); then he worked for more then ten years as a school teacher in remote parts of
Prussia, teaching not only mathematics, but also subjects like botany, calligraphy, and physical education. Finally, almost 40 years old, he became a famous mathematician over night,
and eventually he was appointed as professor at Berlin --- without ever having received a PhD.
\par
This may sound wild, and in some points it simplifies the facts, but it is not far from the truth.
\par
Karl Theodor Wilhelm Weierstra{\ss} was born October 31, 1815, in the village of Ostenfelde in what was then the Prussian province of Westphalia. A street and an elementary school in Ostenfelde are named after
him, and his birth house --- still occupied today --- is listed in a local tourist guide. His father, who worked for for Prussia's customs and taxation authorities, was sent from one post to the next within short
periods of time. For the first fourteen years of Karl Weierstra{\ss}' life, the family was more or less constantly on the move. In 1829, the father obtained an assistant's position at the tax office in the city of
Paderborn (also in Westphalia), and the family could finally settle down. Young Karl enrolled at the local Catholic {\it Gymnasium\/} in Paderborn, where he excelled not only in mathematics, but also in German, Latin, and
Greek. Not only was he an academically strong student --- he was also quite capable of putting his brains to work on much more practical matters: At age 15, he contributed to his family's income by doing the
bookkeeping for a wealthy merchant's widow.
\par
Throughout his life, Weierstra{\ss} senior suffered from the fact that he did not have the right education to rise to a rank in the Prussian civil service that would have really suited his abilities. He had to content himself
with relatively low level positions, not very challenging and not very well paid. Like many a father in this situation, he was determined to prevent such a fate for his bright, eldest son. When Karl graduated in 1834, his father
decided to send him to Bonn to study {\it Kameralistik\/}, i.e.\ a combination of law, finance, and administration. Being a dutiful son, Karl went\ldots
\par
\ldots and did all he could to screw up the life his father had planned for him. He joined a {\it schlagende Verbindung\/}, a particular kind of student fraternity typical of German universities in the nineteenth century.
Besides keeping the brewing industry busy, students in such a fraternity engaged in a peculiar ritual: the {\it Mensur\/}, a swordfight with a particular twist. Unlike in today's athletic competitions, the students fought with
sharp sabers. They wore protective gear which covered most of their bodies --- except the cheeks. During a {\it Mensur\/}, the opponents tried to inflict ghashes on one another's cheeks. The scars resulting from these ghashes
were borne with pride as signs of honor and manhood\footnote{If you find such ideas of honor and manhood absolutely revolting, you're absolutely right.}. Almost two meters tall, quick on his feet, and with strong arms, Karl
Weierstra{\ss} was a fearsome swordsman. His face remained unscarred, and after a while nobody was keen on challenging him anymore. Having escaped from under his father's tutelage, he spent his years at Bonn drinking beer and wielding the saber ---
and seriously studying mathematics. Although he was not enrolled in mathematics, he read some of the most advanced math books of his time. When it was time for him, in 1838, to take his exams, he simply dropped out.
\par
His family was desperate. They had made considerable financial sacrifices to secure a better future for Karl, who had let them down. Having wasted four years of his life, he needed a breadwinning degree --- and that fast.
So, in 1839, he enrolled at the {\it Akademie\/}, a forerunner of today's university, in M\"unster to become a secondary school teacher. Although this was not really a university, but rather a teacher training college, they had one good mathematician 
teaching there: Christoph Gudermann. He is said to have been an abysmal teacher: Very often, he had just one student sitting in his class --- Karl Weierstra{\ss}. In 1840, Weierstra{\ss} graduated. His thesis
was so good that Gudermann believed it to be strong enough for a doctoral degree. However, the {\it Akademie\/} was not {\it really\/} a university: It did not have the right to grant doctorates. So, instead of
receiving a doctorate and starting an academic career, Weierstra{\ss} left the {\it Akademie\/} as a mere school teacher.
\par
His first job (probationary) was in M\"unster. One year later, he was sent to Deutsch Krona\footnote{Now Watcz in Poland.} in the provice of West Prussia as an auxiliary teacher and in 1848 to 
Braunsberg\footnote{Now Braniewn in Poland.} in East Prussia. Of course, he taught mathematics, but also physics, geography, history, German, and --- believe it or not --- calligraphy and physical education.
Besides the demands of working full time as a teacher and having a social life --- remember: he liked beer ---, he found time to do research in mathematics. During his time in Braunsberg, he published
a few papers in his school's yearbook. High school year books are not exactly the place where people are looking for cutting edge research in mathematics, and so, nobody noticed them. Then, in 1854, he published a 
paper entitled ``Zur Theorie der Abelschen Functionen'' in a widely respected journal. I won't even make an attempt to explain what it was about. But unlike his previous work, this one {\it was\/} noticed.
\par
It dawned to mathematicians all over Europe that the man who was probably the leading analyst of his day was rotting in a small East Prussian town, spending most of his time teaching youngsters calligraphy and physical education. 
On March 31, 1854, Weierstra{\ss} finally received a doctorate: an honorary one from the university of K\"onigsberg\footnote{Now Kaliningrad in Russia.}, but nevertheless. In 1856, he accepted a position at the {\it Gewerbeinstitut\/},
an engineering school, and a year later he joined the faculty of the university of Berlin as an adjunct professor. As a teacher, he attracted large audiences. Often, he taught in front of more than 200 students. Finally, in 1860, when
he was almost 50 years old, Weierstra{\ss} was appointed full professor at the university of Berlin. In the years 1873 and 1874, he was {\it Rektor magnificus\/} of the university; in 1875, he became a knight of the
order ``Pour le M\'erite'' in the category of arts and sciences, the highest honor newly unified Germany could bestow upon one of her citizens; and in 1885, on the occasion of his 70th birthday, a commemorative coin was issued.
\par
The years of leading a double life as a secondary school teacher and a mathematical researcher had taken their toll on Weierstra{\ss}' health. A less vigorous man would probably have collapsed under the double burden after a short time only.
From 1850 on, he began to suffer from attacks of dizziness which culminated in a collapse in 1861. He had to pause for a year before he could teach again, and he never recovered fully. In 1890, --- at age 75! --- he retired from teaching
because of his failing health. The last years of his life, he spent in a wheelchair. In 1897, he died.
\par
Weierstra{\ss} published few papers --- he was very critical towards his own work. But although, he was a brilliant researcher, the greatest impact he had on mathematics was as a teacher. At Berlin he repatedly taught a two-year course
on analysis --- the ancestor of all modern introductions to calculus and analysis. Although he never wrote a textbook, notes taken in class by his students have survived and convey an impression of his lectures. The perhaps longest lasting
legacy of these lectures is their ephasis on rigor. When calculus was created in the 17th century, mathematicians did not worry about rigorously proving their results. For example, the 
first derivative $dy/dx$ of a function $y = f(x)$ was thought of as a quotient of two ``infinitesimal'', i.e.\ infintely small, quantities $dy$ and $dx$. Nobody could really tell what infinitely small quantities were supposed to be, but 
mathematicians back then didn't really care. The new mathematics enabled them to solve problems in physics and engineering which had been beyond the reach of human mind before. So why bother with rigor? In the 18th century, mathematicians
went so far as to proclaim that rigor was for philosophers and theologians --- not for mathematicians. But with the lack of rigor, contradictory results cropped up with disturbing frequency, i.e.\ people arrived at formulae that were
obviously wrong. And if a particular formulae determines whether or not a bridge collapses, you don't want it to be wrong. Weierstra{\ss} realized that, if calculus was to rest on solid foundations, its central notion, that of limit, had to
be made rigorous. He introduced the definition which (essentially) is still used today in classes:
\begin{quote} \it
A number $y_0$ is the limit of a function $f(x)$ as $x$ tends to $x_0$ if, for each $\epsilon > 0$, there is $\delta > 0$ such that $|f(x) - y_0 | < \epsilon$ for each $x$ with $| x - x_0 | <\delta$.
\end{quote} 
Students may curse it, but it will not go away.
\par
Weierstra{\ss} was not only an influential lecturer, but also one of the most prolific advisors of PhD theses of all time. There is a database on the internet\footnote{The Mathematics Genealogy Project at
http://hcoonce.math.mankato.msus.edu/.} which lists 31 PhD students of Weierstra{\ss} and 1346 descendants, i.e.\ PhDs of PhDs of PhDs of etc.\ of Weierstra{\ss}. Interestingly, the two of
his former students who generated the most folklore weren't his students in a technical sense.
\par
Sofya Kovalevskaya was a young Russian noble woman who had come to Germany to study mathematics. This alone was not a small feat at a time when the very idea of a woman receiving a university education was
revolutionary. For two years, she studied at Heidelberg, where authorities would not let her enroll officially, but eventually allowed her to attend lectures unofficially (provided the instructor did not object).
Then she moved to Berlin to work with Weierstra{\ss} --- only to find that there she was not even allowed to sit in lectures. This prompted Weierstra{\ss}, by all we know a politically conservative man, to tutor her
privately. Since Kovalevskaya could not receive a doctorate from Berlin, Weierstra{\ss} used his influence to persuade the university of G\"ottingen to award her the degree in 1874. The following nine years, she spent jobhunting. Being a woman 
didn't help. The best job she could get was teaching arithmetic at an elementary school. Finally, in 1883, she was offered a professorship at Stockholm, where she worked until her death in 1891 at the age of 41. Weierstra{\ss}
and Kovalevskaya stayed in touch throughout her mathematical career. After her death, Weierstra{\ss} destroyed their correspondence. This fact, along with Kovalevskaya's striking beauty, gave rise to innuedo that she
may have been more to Weierstra{\ss}, who never married, than just a student. May have been --- but we don't know. This kind of ``folklore'' probably tells more about those who spread it than about those it is supposedly about.
\par
G\"osta Mittag-Leffler also was not strictly speaking Weierstra{\ss}' student. Already being enrolled at the university of Uppsala, Sweden, he came to Berlin in 1875 to attend Weierstra{\ss}' lectures, which turned out to have
an enormous impact on his mathematical development. He then returned to his native Sweden, where he received his doctorate.  Over the years, he became the undisputedly most influential mathematician of his time in Sweden. He made use of 
his clout to overcome the obstacles faced by Sofya Kovalevskaya in connection with her appointment at Stockholm. What Mittag-Leffler is most famous for, however, is not a mathematical accomplishment, but a piece of mathematical folklore. 
To this day, mathematicians suffer quietly from the lack of a Nobel prize in mathematics, and some say, Mittag-Leffler is to blame: According to legend, the first version of Nobel's will mentioned a prize in mathematics. Then, Nobel found 
out that his wife had an affair with Mittag-Leffler. Infuriated by the prospect that his wife's lover could well be the first prize winner, Nobel changed his will and removed the math prize. That's a fine piece of juicy folklore --- but 
nothing more: Like Weierstra{\ss}, Nobel was a lifelong bachelor.
\vfill
\begin{tabbing}
{\it Address\/}: \= Department of Mathematical and Statistical Sciences \\
\> University of Alberta \\
\> Edmonton, Alberta \\
\> Canada T6G 2G1 \\[\medskipamount]
{\it E-mail\/}: \> {\tt vrunde@ualberta.ca} \\[\medskipamount]
{\it URL\/}: \> {\tt http://www.math.ualberta.ca/$^\sim$runde/runde.html}
\end{tabbing} 
\end{document}